\documentclass[11pt]{article}
\usepackage{amscd, amsmath, amssymb, graphicx, floatflt, eepic}

\title{The cone conjecture for Calabi-Yau pairs in dimension two}
\author{Burt Totaro}
\date{  }

\def\Z{\text{\bf Z}}
\def\Q{\text{\bf Q}}
\def\R{\text{\bf R}}
\def\C{\text{\bf C}}
\def\P{\text{\bf P}}

\def\arrow{\rightarrow}

\def\qed{\ QED }

\def\l{\langle}
\def\r{\rangle}

\def\Aut{\text{Aut}}
\def\Bir{\text{Bir}}
\def\PsAut{\text{PsAut}}

\def\Int{\text{Int\,}}
\def\NE{\overline{NE}}
\def\A{\overline{A}}
\def\Pic{\text{Pic}}

\def\l{\langle}
\def\r{\rangle}
\def\H{\overline{H}}
\def\im{\text{im}}
\def\Spec{\text{Spec}}
\def\dim{\text{dim}}
\def\Cox{\text{Cox}}

\begin{document}
\maketitle

\newtheorem{theorem}{Theorem}[section]
\newtheorem{corollary}[theorem]{Corollary}
\newtheorem{lemma}[theorem]{Lemma}
\newtheorem{conjecture}[theorem]{Conjecture}

A central idea of minimal model theory as formulated by Mori
is to study algebraic varieties using convex geometry. The cone
of curves of a projective variety is defined as the
convex cone spanned by the numerical equivalence classes
of algebraic curves;
the dual cone is the cone of nef line bundles. For Fano varieties
(varieties with ample anticanonical bundle),
these cones are rational polyhedral by the cone theorem
\cite[Theorem 3.7]{KM}.
For more general varieties, these cones are not well
understood: they can have infinitely many isolated extremal rays,
or they can be ``round''. Both phenomena occur among
Calabi-Yau varieties such as K3 surfaces, which can
be considered the next simplest varieties after Fano varieties.

The Morrison-Kawamata cone conjecture would
give a clear picture of the nef cone for Calabi-Yau varieties
\cite{Morrison1, Morrison2, Kawamata}. 
The conjecture says that the action of the
automorphism group of the variety on the nef cone 
has a rational polyhedral fundamental domain.
(The conjecture includes an analogous statement about the movable cone;
see section \ref{conjecture} for details.)
Thus, for Calabi-Yau varieties,
the failure of the nef cone to be rational
polyhedral is always explained by an infinite discrete
group of automorphisms
of the variety. It is not clear where these
automorphisms should come from. Nonetheless, the conjecture
has been proved for Calabi-Yau surfaces by Sterk, Looijenga,
and Namikawa \cite{Sterk, Namikawa, Kawamata}, the heart
of the proof being the Torelli theorem for K3 surfaces of
Piatetski-Shapiro and Shafarevich \cite[Theorem 11.1]{BPV}. Kawamata proved
the cone conjecture for all 3-dimensional
Calabi-Yau fiber spaces over a positive-dimensional base
\cite{Kawamata}. The conjecture is wide open for Calabi-Yau 3-folds,
despite significant results by Oguiso and Peternell \cite{Oguiso},
Szendr\" oi \cite{Szendroi}, Uehara \cite{Uehara} and Wilson \cite{Wilson}.

The conjecture was generalized from Calabi-Yau varieties to klt
Calabi-Yau pairs $(X,\Delta)$ in \cite{Totaro}.
Here $\Delta$ is a divisor on $X$ and
``Calabi-Yau'' means that $K_X+\Delta$ (rather than $K_X$)
is numerically trivial. In this paper we prove
the cone conjecture for all klt Calabi-Yau pairs
of dimension 2 (Theorem \ref{rational}), using the geometry
of groups acting on hyperbolic space and 
reduction to the case of K3 surfaces.
This is enough to show that the conjecture is reasonable
in the greater generality of pairs. More concretely,
the theorem gives control over the nef cone
and the automorphism group for a large class of
rational surfaces, including the Fano surfaces as well as many
others. In particular, we get
a good description of when the Cox ring (or total coordinate ring)
is finitely generated in this class of surfaces (Corollary \ref{finite}).

Thanks to Caucher Birkar, Igor Dolgachev, Brian Harbourne,
Artie Prendergast-Smith, and Chenyang Xu for their comments.

\section{The cone conjecture}
\label{conjecture}

In this section, we state the cone conjecture for klt Calabi-Yau pairs
following \cite{Totaro},
and discuss some history and examples.

Varieties are irreducible by definition, and a {\it curve }means
a variety of dimension 1. Our main Theorem \ref{rational} takes the
base field to be the complex numbers, but Conjecture \ref{cone}
makes sense over any field.
For a projective morphism $f:X\arrow S$ of normal varieties
with connected fibers, define $N^1(X/S)$ as the real vector
space spanned by Cartier divisors on $X$ modulo numerical equivalence 
on curves on $X$ mapped to a point in $S$ (that is,
$D_1\equiv D_2$ if $D_1\cdot C=D_2\cdot C$ for all curves
$C$ mapped to a point in $S$). Define a
{\it pseudo-isomorphism }from $X_1$ to $X_2$ over $S$ to be
a birational map $X_1\dashrightarrow X_2$ over $S$ which
is an isomorphism in codimension one.
A {\it small $\Q$-factorial modification }(SQM)
of $X$ over $S$ means a pseudo-isomorphism over $S$ from $X$ to some other
$\Q$-factorial variety with a projective morphism to $S$.
A Cartier divisor $D$ on $X$ is called {\it $f$-nef} (resp.\
{\it $f$-movable, $f$-effective}) if $D\cdot C\geq 0$ for every
curve $C$ on $X$ which is mapped to a point in $S$ (resp.,
if $\text{codim}(\text{supp}(\text{coker}(f^*f_*O_X(D)\arrow O_X(D))))\geq 2$,
if $f_*O_X(D)\neq 0$). 

The canonical divisor is denoted $K_X$. For an $\R$-divisor $\Delta$ on
a normal $\Q$-factorial variety $X$,
the pair $(X,\Delta)$ is {\it klt }if, for all resolutions
$\pi:\widetilde{X}\arrow X$ with a simple normal crossing
$\R$-divisor $\widetilde{\Delta}$
such that $K_{\widetilde{X}}+\widetilde{\Delta}=\pi^*(K_X+\Delta)$,
the coefficients of $\widetilde{\Delta}$
are less than 1 \cite[Definition 2.34]{KM}. It suffices to check
this property on one resolution.
For a complex surface $X$,
the pair $(X,0)$ is klt if and only if
$X$ has only quotient singularities \cite[Proposition 4.18]{KM}.
For later use, we define a pair $(X,\Delta)$ to be {\it terminal }if,
for all resolutions $\pi:\widetilde{X}\arrow X$ as above, 
the coefficients in $\widetilde{\Delta}$ of all exceptional divisors
of $\pi$ are less than 0. (The definition of terminal pairs
puts no restriction on the coefficients of $\Delta$ itself,
although one checks easily (for $\dim(X)\geq 2$)
that they are less than 1;
that is, a terminal pair is klt.) For a surface $X$,
$(X,0)$ is terminal if and only if $X$ is smooth.

 The {\it $f$-nef cone }$\A(X/S)$
(resp.\ the {\it closed $f$-movable cone }$\overline{M}(X/S)$)
is the closed convex cone in $N^1(X/S)$ generated by
the numerical classes of $f$-nef divisors (resp.\ $f$-movable divisors).
The {\it $f$-effective cone }$B^e(X/S)$ is the convex cone,
not necessarily closed, generated by $f$-effective Cartier divisors.
We call $A^e(X/S)=\A(X/S)\cap B^e(X/S)$ and
$M^e(X/S)=\overline{M}(X/S)\cap B^e(X/S)$ the {\it $f$-effective
$f$-nef cone }and the {\it $f$-effective $f$-movable cone, }respectively.
Finally, a {\it rational polyhedral }cone in $N^1(X/S)$ means the closed
convex cone spanned by a finite set of Cartier divisors on $X$.

We say that $(X/S,\Delta)$ is a {\it klt Calabi-Yau pair }if
$(X,\Delta)$ is a $\Q$-factorial klt pair with $\Delta$ effective
such that $K_X+\Delta$ is numerically trivial over $S$. (Our main Theorem
\ref{rational} takes $S$ to be a point.)
Let $\Aut(X/S,\Delta)$ and $\PsAut(X/S,\Delta)$
denote the groups of automorphisms
or pseudo-automorphisms of $X$ over the identity on $S$
that map the divisor $\Delta$ to itself.

\begin{conjecture}
\label{cone}
Let $(X/S,\Delta)$ be a klt Calabi-Yau pair.

(1) The number of $\Aut(X/S,\Delta)$-equivalence classes of faces
of the cone $A^e(X/S)$ corresponding to birational contractions
or fiber space structures is finite. Moreover, there exists
a rational polyhedral cone $\Pi$ which is a fundamental
domain for the action of $\Aut(X/S,\Delta)$ on $A^e(X/S)$ in the sense
that

(a) $A^e(X/S)=\cup_{g\in \Aut(X/S,\Delta)} g_*\Pi$,

(b) $\Int \Pi \cap g_*\Int \Pi=\emptyset$ unless $g_*=1$.

(2) The number of $\PsAut(X/S,\Delta)$-equivalence classes of chambers
$A^e(X'/S)$ in the cone $M^e(X/S)$ corresponding to
marked SQMs $X'\arrow S$ of $X\arrow S$ is finite.
Equivalently, the number of isomorphism classes over $S$
of SQMs of $X$ over $S$ (ignoring the birational
identification with $X$) is finite. Moreover, there exists
a rational polyhedral cone $\Pi'$ which is a fundamental
domain for the action of $\PsAut(X/S,\Delta)$ on $M^e(X/S)$.
\end{conjecture}

For $X$ terminal and $\Delta=0$, Conjecture
\ref{cone} is exactly Kawamata's conjecture
on Calabi-Yau fiber spaces,
generalizing Morrison's conjecture on Calabi-Yau
varieties \cite{Kawamata, Morrison1, Morrison2}.
(The group in part (2) can then be described as $\Bir(X/S)$, since
all birational automorphisms of $X$ over $S$ are
pseudo-automorphisms when $X$ is terminal and
$K_X$ is numerically trivial over $S$.)

Conjecture \ref{cone} implies the analogous statement
for the group of automorphisms or pseudo-automorphisms
of $X$ rather than of $(X,\Delta)$. (That slightly weaker formulation
of the cone conjecture is used in \cite{Totaro}.)

The first statement of part (1) follows from the second statement,
on fundamental domains. Indeed, each contraction of $X$
to a projective variety is given by some semi-ample line
bundle on $X$ (a line bundle for which some positive multiple
is basepoint-free).
The class of such a line bundle in $N^1(X)$ lies in the nef
effective cone.
And two semi-ample line bundles in the interior
of the same face of some cone $\Pi$ determine the same
contraction of $X$, since they have degree zero on the same curves.
Thus the second statement of (1) implies the first. We include
the first statement in the conjecture because one can try to prove it
in some cases where the conjecture on fundamental domains remains
open. The first statements of (1) and of (2) are what Kawamata proves
for Calabi-Yau fiber spaces
of dimension 3 over a positive-dimensional base \cite{Kawamata}.

For $X$ of dimension at most 2, we only need to consider statement (1),
because any pseudo-isomorphism
between normal projective surfaces
is an isomorphism, and every movable divisor on a surface is nef.

Conjecture \ref{cone} would not be true for
Calabi-Yau pairs
that are log-canonical (or dlt) rather than klt. Let $X$
be the blow-up of $\P^2$ at 9 very general points. Let $\Delta$
be the proper transform of the unique smooth cubic curve through
the 9 points; then $K_X+\Delta\equiv 0$, and so $(X,\Delta)$ is
a log-canonical Calabi-Yau pair. The surface $X$ contains
infinitely many $(-1)$-curves by Nagata \cite{NagataRat2},
and so the nef cone is not finite polyhedral. But the
automorphism group $\Aut(X)$ is trivial \cite{Hirschowitz,
Koitabashi} and hence does not
have a finite polyhedral fundamental domain on the nef cone.
There is also an example of a log-canonical Calabi-Yau surface
with rational
singularities (with $\Delta=0$) for which the cone conjecture fails:
contract the divisor $R_1+\cdots+R_4+2R_5\sim -2K_X$ in the surface
$X$ of Dolgachev-Zhang \cite[Example 6.10]{DZ}.

The conjecture also fails if we allow the $\R$-divisor $\Delta$ to have
negative coefficients. Let $Y$ be a K3 surface whose nef cone
is not finite polyhedral, and let $X$ be the blow-up of $Y$ at
a very general point. Let $E$ be the exceptional curve.
Then $(X,-E)$ is klt and $K_X-E\equiv 0$.
The nef cone of $X$
is not finite polyhedral, but $\Aut(X)$ is trivial.

An interesting class of klt Calabi-Yau pairs are the rational
elliptic surfaces (meaning smooth rational surfaces which
are minimal elliptic fibrations over $\P^1$).
The cone conjecture was checked for rational elliptic surfaces
with no multiple fibers and Mordell-Weil rank 8 by
Grassi-Morrison \cite[Theorem 2.3]{GraM}, and for all rational
elliptic surfaces by \cite[Theorem 8.2]{Totaro}. This will
be generalized by our main result, Theorem \ref{rational}. Note that
the cone conjecture for klt pairs (applied to a suitable divisor
$\Delta$ with $K_X+\Delta\equiv 0$) describes the whole nef cone
of a rational elliptic surface $X$ in $N^1(X)$. By contrast,
if we take $\Delta=0$ and apply the cone conjecture
to the elliptic fibration $X\arrow S$,
then we only get information about the relative nef cone in $N^1(X/S)$.
For example, the relative nef cone is (trivially) rational polyhedral
for any rational elliptic surface, whereas the whole nef cone
is rational polyhedral if and only if the Mordell-Weil rank is 0
\cite[Theorem 5.2, Theorem 8.2]{Totaro}.

{\bf Example. }We give an example of a rational surface $Y$,
considered by Zhang \cite[Theorem 4.1]{Zhang} and Blache
\cite[Theorem C(b)(2)]{Blache}, whose nef cone
is a 4-dimensional round cone. The surface $Y$ is klt Calabi-Yau, and so
the cone conjecture is true by our main Theorem \ref{rational}
(known in this case by Oguiso-Sakurai \cite[Corollary 1.9]{OS}). Thus
the automorphism group of $Y$ must be infinite, and in fact it is
a discrete group of isometries of hyperbolic 3-space with quotient
of finite volume. (The quotient of hyperbolic 3-space by an index-24
subgroup of the group here is familiar to topologists as the
complement of the figure eight knot \cite[1.4.3, 4.7.1]{MR}.)

Let $\zeta$ be a primitive cube root of unity.
Let $X$ be the blow-up of $\P^2$ at the 12 points $[1,\zeta^i,\zeta^j],
[1,0,0],[0,1,0],[0,0,1]$ over the complex numbers.
(This is the dual of the famous Hesse configuration
of 9 points lying on 12 lines in the plane \cite[4.6]{DolgachevAbstract}.
Combinatorially, we can identify the Hesse configuration with
the 9 points and 12 lines in the affine plane over $\Z/3$.)
Let $E_1,\ldots,E_9$ be the proper transforms of the 9 lines through
quadruples of
the 12 points; these curves have self-intersection $-3$ in $X$,
and $(X,(1/3)\sum E_i)$ is a klt Calabi-Yau pair. We can contract
the 9 disjoint curves $E_i$ to obtain a klt Calabi-Yau surface $Y$
with 9 singular points of type $(1/3)(1,1)$.
Then $Y$ has Picard number 4, and the point
of this example is that the nef cone of $Y$ is a round cone
(one of the two pieces of $\{x\in N^1(Y): x^2\geq 0\}$). 
That follows by viewing $Y$ as the quotient of an abelian surface
$E\times E$ by $\Z/3$; here $E$ is the elliptic curve $\C/\Z[\zeta]$
and $\Z/3$ acts by $(\zeta,\zeta)$ on $E\times E$. The nef cone
of an abelian surface is always round, and in this case $\Z/3$
acts trivially on $N^1(E\times E)$ so that the nef cone of $Y$
is equal to that of $A$.

The cone conjecture, a theorem in this case, implies that
$\Aut(Y)$ must be a discrete group acting on hyperbolic 3-space
(the ample cone of $Y$ modulo scalars) with finite-volume quotient
(since every finite polytope in hyperbolic space has finite volume,
even if its vertices are at infinity).
In this example, we compute that $\Aut(Y)$ is the group
$(GL(2,\Z[\zeta])/(\Z/3))\ltimes (\Z/3)^2$, and its image
$\Aut^*(Y)$ in $GL(N^1(Y))$ is $PGL(2,\Z[\zeta])$.
A decomposition of hyperbolic 3-space into fundamental domains
for this group looks roughly like the figure
(an analogous picture in the hyperbolic plane).
\begin{floatingfigure}[l]{0.3\textwidth}
\centering
\includegraphics[scale=0.15]{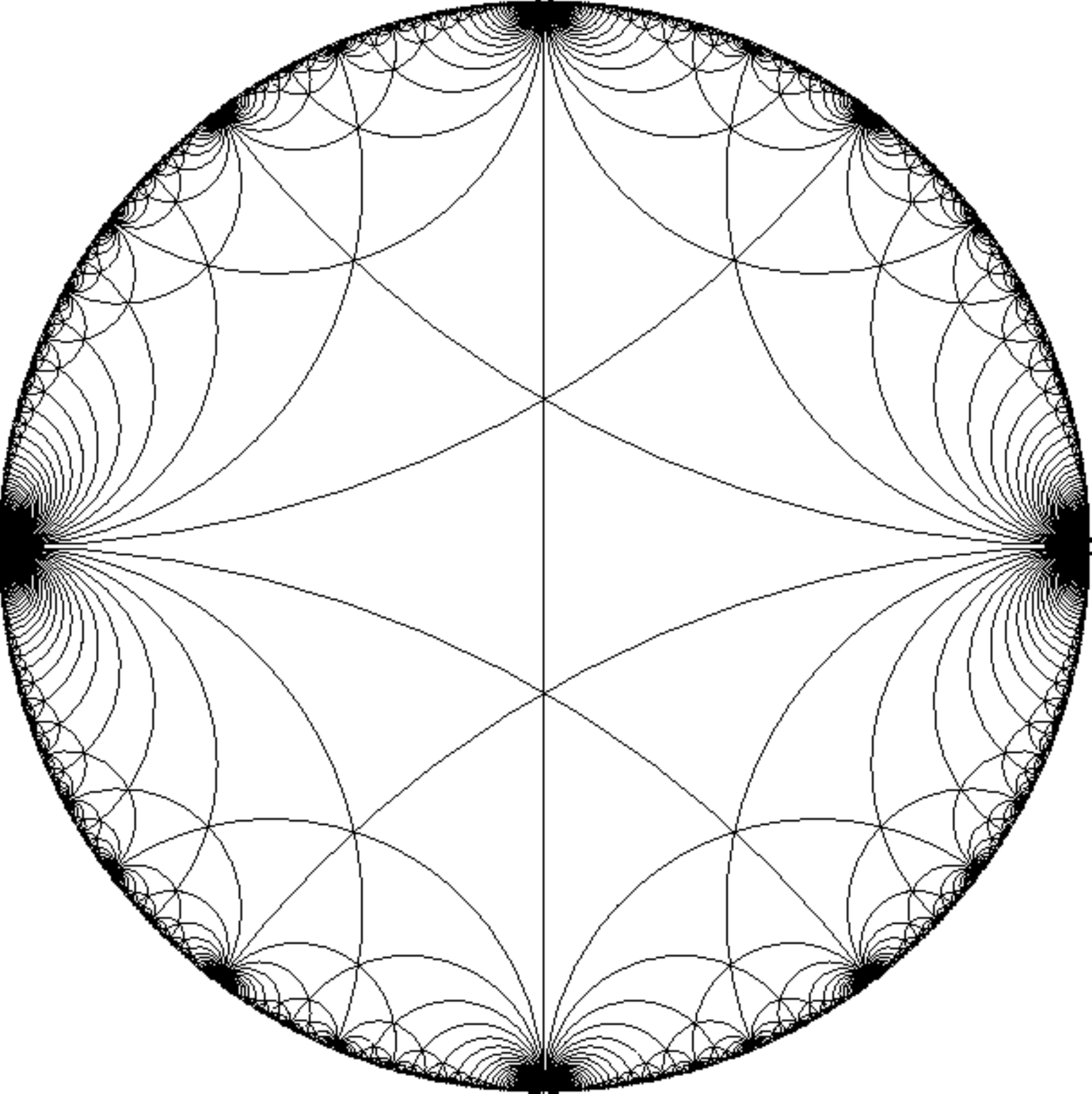}
\end{floatingfigure}
In particular,
any choice of rational polyhedral fundamental domain has a vertex
at the boundary of hyperbolic space as in the figure,
because there are rational points on the boundary of the nef cone
(corresponding
to elliptic fibrations of $Y$). The automorphism group of the
smooth rational surface $X$ is the same as that of $Y$,
but now acting on hyperbolic space of dimension $\rho(X)-1=12$.
The nef cone of $X$ modulo scalars is not all of hyperbolic space, because $X$
contains curves of negative self-intersection (including
infinitely many $(-1)$-curves). The cone
conjecture for $X$ is not immediate from that for $Y$, but it also
follows from Theorem \ref{rational}.

Another feature of this example is that
$PGL(2,\Z[\zeta])$ is the unique non-cocompact group
of orientation-preserving isometries of hyperbolic 3-space
of minimum covolume (about 0.085) \cite{Meyerhoff}.

\section{Hyperbolic geometry}
\label{hyperbolic}

For $X$ of dimension 2, we can view the nef cone
of $X$ modulo scalars as a convex subset of hyperbolic space.
Using hyperbolic geometry, we show that a
weak statement in the direction of the cone conjecture automatically
implies the full conjecture (Lemma \ref{dirichlet}).

Throughout this section, let $X$ be a klt surface.
Let $G$ be any subgroup of the image $\Aut^*(X)$ of $\Aut(X)$
in the orthogonal
group $O(S)$, where $S:=\im(\Pic(X)\arrow N^1(X))$ is a lattice
with signature $(1,n)$ by the Hodge index theorem. More precisely,
$G$ is contained in the index-two subgroup $O^+(S)$ of isometries
that preserve the positive cone $\{x\in S_{\R}: x^2>0,\text{ }
A\cdot x>0\}$ (where $A$ is any ample line bundle on $X$).

It is convenient to identify the quotient of the positive
cone modulo positive real scalars with hyperbolic space $H$
of dimension $n$.
This is one of the standard models for hyperbolic space.
Explicitly, the distance between two points in hyperbolic space,
represented by vectors $x$ and $y$ in the positive cone
with $\l x,x\r=\l y,y\r=1$,
is given by $\cosh d(x,y)=\l x,y\r$ \cite[p.~26]{Vinberg}.
The closure of the positive cone minus the origin, modulo scalars,
is the compactification $\H$, which has a conformal model
as a closed ball. 
We define the rational points $\partial H(\Q)$
in the boundary $\partial H=\H-H$ to be the images
of nonzero rational points in the boundary of the positive cone.

\begin{floatingfigure}[l]{0.3\textwidth}
\centering
\includegraphics[scale=0.2]{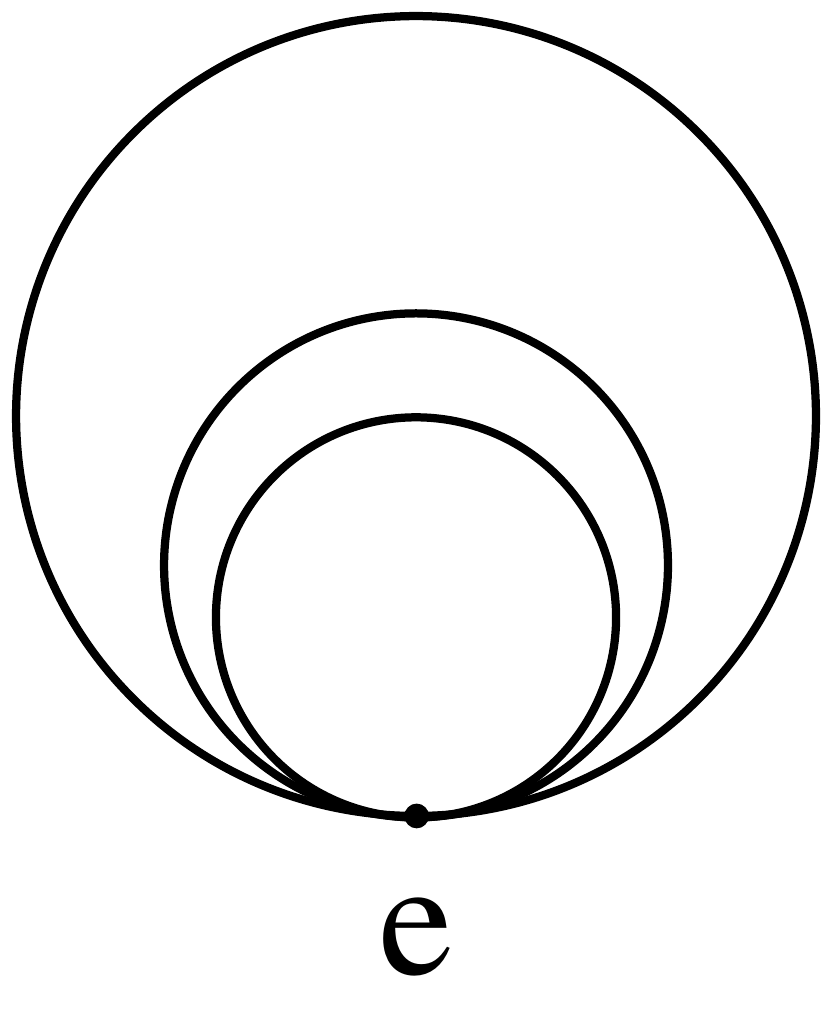}
\end{floatingfigure}
The {\it horoballs }at a point in $\partial H$ represented
by a vector $e$ with $\l e,e\r =0$ are defined as the images
in $H$ of the cones $\l x,e\r\leq a\sqrt{\l x,x\r}$ in the positive cone,
for positive constants $a$. In the picture of $\H$ as a ball in $\R^n$,
these are the balls inside $\H$ tangent to the boundary at $e$.
Given one horoball $U_e=\{x\in \H: \l x,e\r \leq a\sqrt{\l x,x\r}\}$
and a positive number $c$,
we write $cU_e$ for the horoball $cU_e=
\{ x\in\H: \l x,e\r \leq ca\sqrt{\l x,x\r}\}$. Every point of the
horosphere $\partial(cU_e)-\{e\}$ is at
distance $\log(c)$ from $\partial(U_e)-\{e\}$ in hyperbolic space.

For each rational point in $\partial H$, the image of a primitive vector
$e$ in $S$ with $\langle e,e\rangle =0$, the stabilizer of $e$ in
the orthogonal group $O^+(S)$ contains as a subgroup
of finite index the group $e^{\perp}/\Z e\cong \Z^{n-1}$
of ``strictly parabolic''
transformations at $e$. Explicitly, for $x \in e^{\perp}$, these
transformations can be defined as
$$\alpha_x(y)=y+\langle y,e\rangle x-\big[ \langle x,y\rangle
+\frac{1}{2}\langle x,x\rangle \langle y,e\rangle \big] e.$$
Geometrically, this group preserves each of the horoballs at
$e\in \partial H$.
The boundary of each horoball inside hyperbolic $n$-space $H$
is isometric to Euclidean space $\R^{n-1}$ \cite[p.~87]{Vinberg},
and this group
$\Z^{n-1}\subset O^+(S)$ acts as a discrete cocompact group
of translations on $\R^{n-1}$.

We use the following elementary inequality.

\begin{lemma}
\label{ineq}
Let $e_1,e_2$ be points in the boundary of the positive cone.
Let $x$ be a point in the interior of the positive cone,
scaled to have $\l x,x\r=1$. Then $\l e_1,e_2\r\leq 2
\l x,e_1\r \l x,e_2\r$.
\end{lemma}

{\bf Proof. }Since the intersection form on the lattice $S$
has signature $(1,n)$, the intersection form on $x,e_1,e_2$
has signature $(1,2)$ if it is nondegenerate. So the determinant
of the matrix of inner products, 
$\begin{pmatrix} 1 & \l x,e_1 \r & \l x,e_2\r \\
\l x,e_1\r & 0 & \l e_1,e_2\r \\
\l x,e_2\r & \l e_1,e_2\r & 0
\end{pmatrix}$,
is nonnegative. This gives the inequality we want. \qed

A general fact about discrete groups of isometries of hyperbolic
space is that there is a set of ``precisely invariant'' horoballs
for the group action \cite[p.~171]{Vinberg}; this is used to
compactify quotients of hyperbolic space by adding ``cusps'' or
boundary components.
For the orthogonal group $O^+(S)$, here is
an elementary proof of what we need.
We know that the bilinear form on the lattice $S$ takes
values in the integers.
For each rational
point in $\partial H$, let $e$ be the unique primitive vector
in $S$ representing this point. Let $U_e$ be the horoball
at $e$ defined as
$$U_e:=\{x\in H: \langle x,e\rangle \leq \frac{1}{2}\sqrt{
\langle x,x\rangle } \}.$$
(Any constant less than $1/\sqrt{2}$ could be used in place of
$1/2$.)
For any two distinct primitive
vectors $e_1$ and $e_2$ in the boundary of the positive cone,
we have $\langle e_1,e_2\rangle >0$ and hence $\langle e_1,e_2
\rangle \geq 1$. By Lemma \ref{ineq}, it follows that
$U_{e_1}\cap U_{e_2}=\emptyset$. Thus, we have attached
a horoball to each rational point of $\partial H$ in such a way
that $O^+(S)$ permutes these horoballs (clearly), distinct
horoballs are disjoint, and the subgroup of $O^+(S)$ mapping
a given horoball $U_e$ into itself is exactly the stabilizer
of $e$ in $O^+(S)$.

Let $y$ be a rational point in the interior of the nef cone whose
stabilizer group in $G$ is trivial. We define the associated
{\it Dirichlet domain }$D$ as the closed convex cone
$$D=\{ x\in \A(X): \langle x,y\rangle \leq \langle x,gy\rangle
\text{ for all }g\in G\}.$$
Modulo scalars, we can view $D$ as a closed convex subset of $\H$;
its intersection with hyperbolic space is the set of points
in our convex set $\A(X)$ that are at least as close to $y\in H$,
in the hyperbolic metric,
as they are to any other point in the $G$-orbit of $y$.

\begin{lemma}
\label{dirichlet}
Suppose we are given a $G$-invariant collection
of rational polyhedral cones
whose union is $A^e(X)$. Suppose that 
the cones fall into only finitely
many $G$-orbits. Then the Dirichlet domain $D\subset \A(X)$ associated
to a rational point with trivial stabilizer
in the interior of one of these cones 
is rational polyhedral, it is contained in $A^e(X)$,
and $A^e(X)=\cup_{g\in G}\: gD$.
\end{lemma}

It is clear from the definition of a Dirichlet domain that
the interiors of $D$ and $gD$ are disjoint, for $g\neq 1$
in $G$. So the conclusion of this lemma says that the cone conjecture
holds, starting from a weaker assumption (where we have
finitely many $G$-orbits of cones instead of just one,
and no requirement about how they intersect).
This implication is easy to believe, but the proof seems to require
a fair amount of hyperbolic geometry.
For comparison, there are some geometrically finite discrete groups
acting on hyperbolic space for which most Dirichlet domains
have infinitely many faces \cite[p.~173]{Vinberg}.

\begin{corollary}
\label{locfin}
Under the assumption of Lemma \ref{dirichlet} (in particular,
if the cone conjecture holds), the collection of rational polyhedral
cones $gD$ that covers $A^e(X)$ is locally finite in the interior
of the positive cone. 
\end{corollary}

{\bf Proof of Corollary \ref{locfin}. }Because $G\subset O^+(S)$,
$G$ acts discretely on hyperbolic space. So
the $G$-orbit of the point $y$ defining the Dirichlet
domain is discrete in hyperbolic space. As a result,
the decomposition $A^e(X)=\cup_{g\in G}\: gD$ of Lemma \ref{dirichlet}
is locally finite in hyperbolic space, or equivalently in the positive
cone. \qed

{\bf Proof of Lemma \ref{dirichlet}. }We work modulo scalars,
in hyperbolic space
$H$ and its compactification $\H$. We first show that
$A^e(X)$ is contained in $\cup_{g\in G}\: gD$. Since $G$ acts discretely
on hyperbolic space, the orbit $Gy$ is a discrete subset of $H$. 
So every point in $H$ has at least one closest point in $Gy$,
and thus $A^e(X)\cap H$ is contained in $\cup_{g\in G}\:  gD$
by definition of $D$. It remains to consider a point $p$
in $A^e(X)\cap \partial H$.

To show that $p$ is in $gD$ for some $g\in D$ means to show that
the horosphere at $p$ through $gy$ contains no $G$-translate of $y$
in its interior. I claim that the
intersection of the $G$-orbit of $y$ with any given horoball
at $p$ is a finite union of $G_p$-orbits. Since $G_p$ preserves
each horoball at $p$, this claim will clearly imply
that there is some point $gy$ ``closest'' to $p$ (i.e. on the smallest
horosphere at $p$) and hence that $p$ is in $\cup_{g\in G}\: gD$,
as we want.

This would be clear if $Gy$ meets the ``precisely
invariant'' horoball $U_p$ we constructed (then $Gy\cap U_p$
would be a single $G_p$-orbit), but it might not. For any $c>1$
and any $q\in \partial H(\Q)$,
write $cU_q$ for the bigger horoball at distance $\log(c)$ from $U_q$.
These bigger horoballs are no longer disjoint, but I claim
that for any $c>1$, only finitely many of the horoballs
$cU_q$ contain the point $y$ in $H$. Indeed, it is equivalent
to show that the ball of radius $\log(c)$ around $y$ meets
only finitely many of the horoballs $U_p$, which is clear.

We can now prove the earlier claim that $Gy$ meets any given
horoball $cU_p$ in finitely many $G_p$-orbits. Enlarging $c$
if necessary, we can assume
that $y$ itself is in $cU_p$. To say that $gy$ is in $cU_p$
means that $y$ is in $cU_{g^{-1}p}$. We have shown that there
are only finitely many points $q$ in $\partial H(\Q)$ with
$y$ in $cU_q$, so the elements $g\in G$ with
$gy$ in $cU_p$ fall into finitely many $G_p$-orbits,
as we want. As discussed earlier, it follows that
$A^e(X)$ is contained in $\cup_{g\in G}\: gD$.

We now begin to prove that $D$ is rational polyhedral.
It is clear that $D$ is locally rational polyhedral
inside hyperbolic space $H$. Next, 
define a {\it chimney }at a point $p$ in $\partial H$ 
to be the convex hull of
$p$ together with some bounded convex polytope in $\partial V_p-\{p\}
\cong \R^{n-1}$, for some horoball $V_p$ at $p$.
We will show that
for each point $p$ in $D\cap \partial H(\Q)$,
$D\cap U_p$ is the chimney over a bounded rational
polyhedron in $\partial U_p-\{p\}\cong \R^{n-1}$.
The assumption that $p$ is in $D$
means that $y$ is on the closest horosphere $\partial (cU_p)$ to $p$
among all the points in the $G$-orbit of $y$.

We showed two paragraphs back that the $G$-orbit
of $y$ meets the closed horoball $cU_p$ in finitely
many $G_p$-orbits. The strictly parabolic elements of $G_p$
form an abelian subgroup $T_p$ of finite index, and so
the $G$-orbit of $y$ meets $cU_p$ in finitely
many $T_p$-orbits. By our current assumptions, these $T_p$-orbits are all
on the horosphere $\partial (cU_p)$. Here $T_p\cong \Z^a$ acts
discretely by translations on $\partial (cU_p)-\{p\}\cong \R^{n-1}$,
for some $0\leq a \leq n-1$.

The assumption of this lemma gives
a collection of polytopes in $\H$ whose union
is $A^e(X)$. I claim that the polytopes in the given collection that meet
$U_p$ fall into only finitely many $T_p$-orbits. Since the collection
contains only finitely many $G$-orbits of polytopes,
it suffices to consider the $G$-translates $gP$ of a single
polytope $P$ in the collection. We know that $P$ meets
only finitely many of the horoballs $U_x$ for $x\in \partial H(\Q)$,
because $P$ minus the horoballs $U_x$ for $x$ in $P\cap \partial H$
is bounded in $H$. For $P$ to meet a horoball $U_{g^{-1}(p)}$
means exactly that the translate $gP$ meets $U_p$. Since
the stabilizer of $U_p$ is equal to $G_p$, 
the translates $gP$ that meet $U_p$ fall into finitely many
$G_p$-cosets, hence finitely many $T_p$-cosets as we want.

After shrinking $U_p$,
we can assume that all the polytopes $P$
that meet $U_p$ have $p$ as a vertex, and that
$P\cap U_p$ is the convex hull of $p$ and a bounded
rational polytope in $\partial U_p-\{p\}=\R^{n-1}$.
Since $A^e(X)$ is the union of the given collection of polytopes, we conclude
that $A^e(X)\cap U_p$ is the convex hull of $p$ and finitely
many $T_p$-orbits of bounded rational polytopes in $\partial
U_p-\{p\}=\R^{n-1}$.

Let $D_p$ be the Dirichlet domain in $A^e(X)\cap U_p$
associated to $y$ and the other points $Gy\cap \partial
(cU_p)$, that is, 
$$D_p=\{x\in A^e(X)\cap U_p: \l x,y\r \leq \l x,gy\r
\text{ for all }g\in G\text{ with }gy\in \partial(cU_p)\}.$$
Then $D\cap U_p\subset D_p$. We know that $Gy\cap \partial
(cU_p)$ is a finite union of $T_p$-orbits, while $A^e(X)\cap U_p$
is the convex hull of $p$ and finitely many $T_p$-orbits of 
bounded rational polytopes in $\partial U_p-\{p\}=\R^{n-1}$.
It is elementary, then, that $D_p$ is the convex hull
of $p$ and a bounded rational polytope in
$\partial U_p-\{p\}=\R^{n-1}$. Since $D$ is convex
and $D\cap U_p=D_p$, it follows that $p$ is an isolated
point of $D\subset \partial H$.

We showed earlier in this proof that the intersection of the
$G$-orbit of $y$ with any horoball $dU_p$ contains only
finitely many $T_p$-orbits. Therefore there is a $d>c$
such that $Gy\cap dU_p=Gy\cap \partial(cU_p)$. Using the following
geometric Lemma \ref{hyperplane},
it follows that there is a small $a>0$ such that
$D\cap aU_p=D_p\cap aU_p$. Therefore the Dirichlet domain
$D$ is rational polyhedral near each point of $D\cap
\partial H(\Q)$.

\begin{lemma}
\label{hyperplane}
Let $p$ be a point in $\partial H(\Q)$. Let $d>c$ be positive real numbers.
Let $y$ be a point
in $\partial(cU_p)-\{p\}=\R^{n-1}$  and let $D_p$ be the convex hull of $p$
and a bounded neighborhood of $y$ in $\partial(cU_p)-\{p\}$.
Then there is a small $a>0$ such that for every point $z$
in hyperbolic space outside $dU_p$, the hyperplane
of points equidistant from $y$ and $z$ does not meet
$D_p\cap aU_p$.
\end{lemma}

\begin{floatingfigure}[l]{0.3\textwidth}
\centering
\includegraphics[scale=0.3]{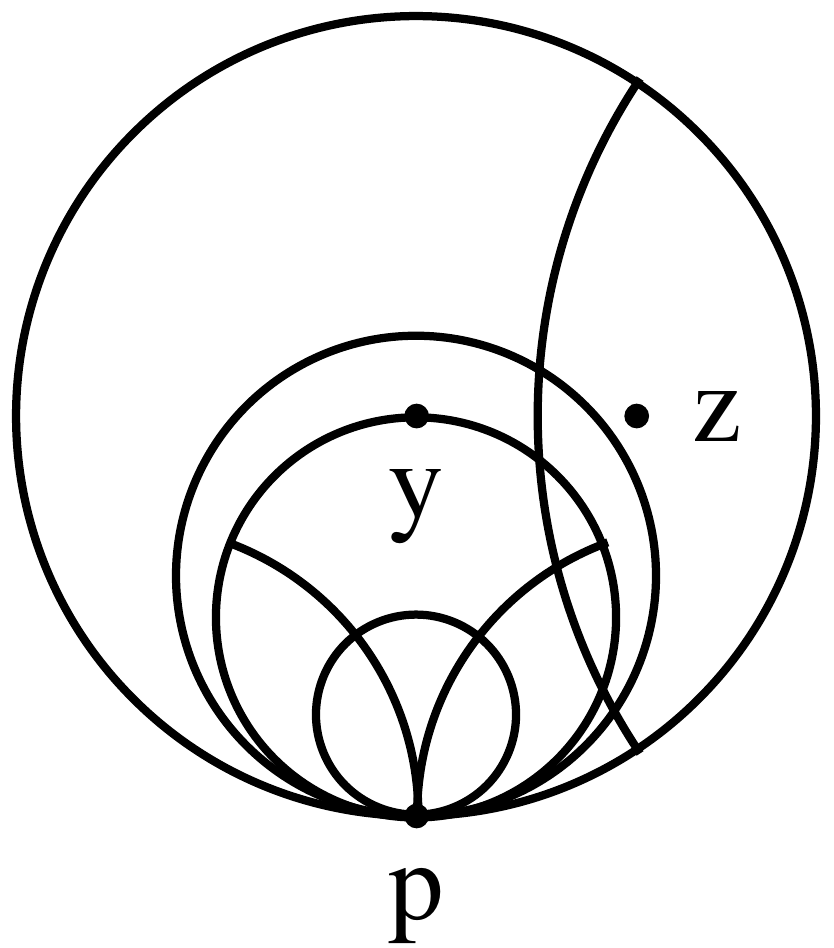}
 \end{floatingfigure}
{\bf Proof. }Suppose this is false. Then there is a sequence
$a_i\arrow 0$ and points $z_i$ in $H-dU_p$ such that the
hyperplane bisecting $y$ and $z_i$ meets $D_p\cap a_iU_p$.
We know that $y$ and $p$ are on the same side of this hyperplane,
because $y$ is on a closer horoball to $p$ than $z_i$ is.
I claim that this hyperplane meets $D_p\cap \partial(cU_p)$ (which
contains $y$); if not, then $D_p\cap \partial(cU_p)$ and $p$
would all be on the same side of the hyperplane, and so
their convex hull $D_p$ would be on that side, contradicting
that the hyperplane meets $D_p\cap a_iU_p$.

Since the hyperplane bisecting $y$ and $z_i$ meets
$D_p\cap a_iU_p$ with $a_i\arrow 0$ and also the bounded
set $D_p\cap \partial(cU_p)$ in hyperbolic space, a subsequence of these
hyperplanes converges to a hyperplane $A$ through $p$ 
that meets $D_p\cap \partial(cU_p)$. So the points $z_i$ must
converge to the point $z$ in hyperbolic space such that $A$
bisects $y$ and $z$. Since $A$ goes through $p$, $z$ lies
on the same horosphere $\partial(cU_p)$ that $y$ does. This
contradicts that all the points $z_i$ are outside the
horoball $dU_p$, where $d>c$. \qed

We return to the proof of Lemma \ref{dirichlet}.
We have shown that the Dirichlet domain $D$ is locally
rational polyhedral inside hyperbolic space, and also
near each point of $D\cap \partial H(\Q)$. To say that $D$ is a finite
rational polytope, it remains to show that $D$ contains no irrational
point of $\partial H$. (If we had defined $D$ as a subset of
$A^e(X)$ instead of $\A(X)$, we would have the equivalent
problem of showing that $D$ is closed in $\H$.)

We can at least say that for each polytope $P$ in the given
collection, $P\cap D$ is a finite rational polytope, because
the vertices of $P$, if in $\partial H$, are rational. Since
the polytopes $P$ are permuted by $G$, $P\cap gD$ is 
a finite rational polytope for each $g\in G$. Moreover,
I claim that $P$ is covered by finitely many of the Dirichlet
domains $gD$. This is clear on a bounded subset of hyperbolic
space, and so it suffices to show that a neighborhood of each
point $p\in P\cap \partial H$ is covered by finitely many
$G$-translates of $D$. This is clear from our description
of $A^e(X)\cap aU_p$, for $a>0$ small: it is contained
in the union of the Dirichlet domains associated to finitely many
$T_p$-orbits of points in the $G$-orbit of $y$, where
each of these domains is the convex hull of $p$ and a bounded
subset of $\partial(aU_p)-\{p\}=\R^{n-1}$.
Since $P\cap \partial(aU_p)-\{p\}$
is bounded, it is contained in the union of finitely many
of these Dirichlet domains.
Thus we have shown that each polytope $P$ in the given
collection is covered by finitely many translates $gD$.

We know that there is a finite union $Q$ of polytopes in our
collection such that $A^e(X)=\cup_{h\in G}\: hQ$. By the previous
paragraph, there is a finite subset $S$ of $G$ such that
$Q=\cup_{g\in S}\: (Q\cap gD)=\cup_{g\in S}\: g(g^{-1}Q\cap D).$
Therefore $A^e(X)=\cup_{h\in G}\: h(\cup_{g\in S}\: g^{-1}Q\cap D)$.
But a point in the interior of $D$ is not in the $G$-orbit
of any other point of $D$, by definition of the Dirichlet
domain. So this equality, applied to points in $\Int(D)$,
implies that $\cup_{g\in S}\: g^{-1}Q\cap D$ contains
$\Int(D)$. But this is a finite union of finite rational polytopes
in $\H$, and so it is closed in $\H$ and hence contains
the closed convex set $D$. Thus the Dirichlet domain
$D\subset \A(X)$ is a rational polyhedral cone, as want.
Finally, since $D$ is equal to $\cup_{g\in S}\: g^{-1}Q\cap D$
where the cones $g^{-1}Q$ are contained in $A^e(X)$,
$D$ is contained in $A^e(X)$. \qed

\section{Klt Calabi-Yau surfaces}

In this section, we prove the cone conjecture for klt Calabi-Yau surfaces
(as opposed to pairs), Theorem \ref{surface}.
This result was stated by Suzuki \cite{Suzuki}. Suzuki's ideas
were inspiring for this paper, but the proof
there is incomplete.

To see the difficulty, let $Y$ be a K3 surface with a node.
Then there are two possible types of curves on $Y$
with negative self-intersection, those with self-intersection $-2$
disjoint from the node and those with self-intersection $-3/2$.
(If $X\arrow Y$ is the minimal resolution with exceptional curve $E$,
the second type of curve is the image of a $(-2)$-curve on $X$ that
meets $E$ transversely in one point.)
The second type is missing in \cite{Suzuki} (see definition of the
set $N'$ and the reflection group $\Gamma$). This makes a difference
because the reflection in a $(-3/2)$-curve does not preserve
the $\Z$-lattice $S=\{x\in \Pic(X): x\cdot E=0\}$. Worse, the angle
between two such reflections need not be a rational multiple
of $\pi$ (take $(-3/2)$-curves $C_1$ and $C_2$ through the same node
of $Y$ with $C_1\cdot C_2=1/2$; this is what happens if the proper
transforms of $C_1$ and $C_2$ are disjoint on the minimal resolution $X$).
So the group
generated by reflections in $(-3/2)$-curves need not be discrete
in $GL(S_{\R})$. As a result, the nef cone of $Y$ need not be a
Weyl chamber for any reflection group acting on the positive
cone. So Lemma 2.4 and Proposition 2.5 in \cite{Suzuki}
do not work.

Our proof of Theorem \ref{surface}, applied to the example
of a K3 surface with a node, works instead by reducing to the
minimal resolution, using hyperbolic geometry in the form of
Lemma \ref{dirichlet}. 
We first note the following consequence
of the abundance theorem in dimension 2.

\begin{lemma}
\label{abun}
Let $(X,\Delta)$ be a klt Calabi-Yau pair of dimension 2.
Then any nef effective $\R$-divisor on $X$ is semi-ample.
\end{lemma}

{\bf Proof. }Since $(X,\Delta)$ is a klt pair
and $K_X+\Delta$ is nef (being numerically trivial),
$K_X+\Delta$ is semi-ample by the abundance theorem
in dimension 2 \cite{Fujita, FM}.
Therefore $K_X+\Delta$ is $\R$-linearly
equivalent to zero. Next,
for any nef effective $\R$-divisor $D$,
$(X,\Delta+\epsilon D)$ is a klt pair for $\epsilon>0$ small,
and $K_X+\Delta+\epsilon D$ is nef. By abundance again,
$K_X+\Delta+\epsilon D$
is semi-ample. That is, $\epsilon D$ is semi-ample. \qed

\begin{theorem}
\label{surface}
The cone conjecture holds for any klt Calabi-Yau surface.
\end{theorem}

{\bf Proof. }Let $Y$ be a klt Calabi-Yau surface. Let
$I=I(Y)$ be the global index of $Y$, that is, the least positive
integer such that $IK_Y$ is Cartier and linearly equivalent to zero,
and let $Z=\Spec(\oplus_{i=0}^{I-1}O_Y(-iK_Y))\arrow Y$
be the global index-one cover of $Y$. Then $Z$ is a
surface with Du Val singularities that has trivial
canonical bundle, and $Y$
is the quotient of $Z$ by an action of $\Z/I$. Let $M$ be the minimal
resolution of $Z$. The smooth surface $M$ has trivial
canonical bundle and hence is a K3 surface or abelian surface. By uniqueness
of the minimal resolution, $\Z/I$ acts on $M$; let $X$ be the quotient
surface.
$$\begin{CD}
M @>>\Z/I> X \\
@VVV @VVV \\
Z @>>\Z/I> Y
\end{CD}$$

By Sterk-Looijenga-Namikawa, we know the cone conjecture
for the smooth Calabi-Yau surface $M$ \cite{Sterk},
\cite[Theorem 2.1]{Kawamata}.
Oguiso-Sakurai proved an analogous statement
for finite groups acting on smooth Calabi-Yau surfaces,
which gives the cone conjecture for $X=M/(\Z/I)$
\cite[Corollary 1.9]{OS}.
That is, there is a rational polyhedral cone $B\subset
A^e(X)$ which is a fundamental domain for the action
of $\Aut(X)$ on the nef effective cone $A^e(X)$.
The theorem follows from Lemma \ref{resolution},
where we take $\Delta=0$. \qed

\begin{lemma}
\label{resolution}
Let $X\arrow Y$ be a proper birational morphism of klt
surfaces. Let $\Delta$ be an $\R$-divisor on $X$ and $\Delta_Y$
its pushforward to $Y$.
If $\Aut(X,\Delta)$ has a rational polyhedral fundamental domain
on the nef effective cone of $X$,
then $\Aut(Y,\Delta_Y)$ has a rational polyhedral
fundamental domain on the nef effective cone of $Y$.
\end{lemma}

{\bf Proof. }The cone of curves $NE(X)$ is defined as the 
convex cone spanned by the classes of curves in $N_1(X)=N^1(X)^*$.
Let $F_0$ be the face of $\NE(X)$
spanned by the curves in $X$ that map to a point in $Y$.
Then the nef cone $\A(X)$ has nonnegative pairing
with $F_0$, and the nef cone of $Y$ is $\A(Y)=\A(X)\cap F_0^{\perp}$;
thus $\A(Y)$ is a face of $\A(X)$. Likewise, the nef effective
cone of $Y$ is $A^e(Y)=A^e(X)\cap F_0^{\perp}$, as one immediately
checks. (In one direction, the image in $Y$ of an effective
divisor on $X$ is effective; in the other, the pullback to $X$
of an effective $\Q$-divisor on $Y$ is effective.)

The subgroup $H$
of $G=\Aut(X,\Delta)$
that maps the face $F_0$ of curves contracted by $X\arrow Y$
into itself is a subgroup of $\Aut(Y,\Delta_Y)$. Equivalently, $H$
is the subgroup of $G$
that maps the face $\A(Y)$ of $\A(X)$ into itself. If we prove
the cone conjecture for this subgroup of $\Aut(Y,\Delta_Y)$,
the statement for the whole group $\Aut(Y,\Delta_Y)$ follows.

We know that there is a rational polyhedral cone
$B$ for $G$ acting on $A^e(X)$, and so 
$A^e(X)=\cup_{g\in G}\: gB$.
It follows that $A^e(Y)=\cup_{g \in G}\: gB\cap
F_0^{\perp}$. Here each set $gB\cap F_0^{\perp}$ is a
rational polyhedral cone contained in $A^e(Y)$.

We will show that these cones fall into finitely many
orbits under $H\subset \Aut(Y,\Delta_Y)$. 
For an element $g$ of $G$, $gB\cap F_0^{\perp}$
is a face of $gB$ (possibly just $\{0\}$). So we can divide the
nonzero intersections $gB\cap F_0^{\perp}$ into finitely
many classes corresponding to the faces $B_i$ of $B$
such that $gB\cap F_0^{\perp}=gB_i$. If the face of $\A(X)$
spanned by $gB_i$ is all of $\A(Y)$, then any automorphism
of $(X,\Delta)$ which maps one such face to another preserves $\A(Y)$
and hence is an automorphism of $Y$, so we have the desired
finiteness.

In general, the face $\A(X)$ spanned by $gB_i$ might be a face
of $\A(Y)$, not all of it. Assume that $gB_i\neq \{0\}$.
Consider the contraction $X\arrow Z$ given
by a $\Q$-divisor in the interior of $gB_i$; this makes sense
because $gB_i$ is contained in the nef effective cone of $X$
and every nef effective $\Q$-divisor on $X$ is semi-ample
by Lemma \ref{abun}. Then we have
$X\arrow Y\arrow Z$ and $Z$ is not a point because
$gB_i\neq \{0\}$. Since $X\arrow Z$ has fiber dimension at most 1,
there are only finitely many numerical equivalence classes
of curves in $X$ contracted by $X\arrow Z$. Therefore there
are only finitely many intermediate factorizations $X\arrow Y'\arrow
Z$. So the automorphisms of $(X,\Delta)$ that map the face $gB_i$ of $\A(X)$
into itself (and hence preserve the contraction $X\arrow Z$) send the face
$\A(Y)$ into only finitely many other faces $\A(Y')$. As a result,
among all $g\in G$ such that $gB\cap F_0^{\perp}=gB_i$,
the cones $gB_i$ fall into only finitely many orbits under
$H$. (We assumed $gB_i\neq \{0\}$, but this conclusion
is also true when $gB_i=\{0\}$.)
Thus $A^e(Y)$ is the union of the rational
polyhedral cones $gB\cap F_0^{\perp}$, and these cones
fall into finitely many orbits under $H\subset \Aut(Y,\Delta_Y)$.

By Lemma \ref{dirichlet}, these properties imply
the cone conjecture for $(Y,\Delta_Y)$. \qed

\section{Klt Calabi-Yau pairs of dimension 2}

\begin{theorem}
\label{rational}
Let $(X,\Delta)$ be a klt Calabi-Yau pair of dimension 2
over the complex numbers. Then Conjecture \ref{cone}
is true. That is, the action of $\Aut(X,\Delta)$ on the nef effective cone
has a rational polyhedral fundamental
domain. As a result, the number of $\Aut(X,\Delta)$-equivalence classes 
of faces of the nef effective cone corresponding to birational
contractions or fiber space structures is finite.
\end{theorem}

We remark that Harbourne's ``K3-like rational surfaces''
have similar finiteness properties \cite{Harbourne}, although
they never have a divisor $\Delta$ with $(X,\Delta)$ klt Calabi-Yau.

{\bf Proof of Theorem \ref{rational}.  }Let
$(\widetilde{X},\widetilde{\Delta})$ be the terminal model
of $(X,\Delta)$. That is,
we have a birational projective morphism $\pi:\widetilde{X}\arrow X$,
$\widetilde{\Delta}$ is effective, $K_{\widetilde{X}}+\widetilde{\Delta}
=\pi^*(K_X+\Delta)$, and $(\widetilde{X},\widetilde{\Delta})$ is terminal.
(Informally, the terminal model of $(X,\Delta)$
is the maximum blow-up of $X$
such that $\widetilde{\Delta}$
is effective \cite[Corollary 1.4.3]{BCHM}.) Thus
$(\widetilde{X},\widetilde{\Delta})$ is a terminal Calabi-Yau pair,
and in particular $\widetilde{X}$ is smooth.
The cone conjecture
for $(\widetilde{X},\widetilde{\Delta})$
implies it for $(X,\Delta)$, by Lemma \ref{resolution}. So we can assume
that $X$ is smooth and $(X,\Delta)$ is a terminal Calabi-Yau pair.

{\bf Example. }The terminal model of a pair
$(X,0)$ of dimension 2 may involve more blowing up
than the usual minimal resolution of $X$.
For a smooth surface $X$ and a divisor
$\Delta$ consisting of two smooth curves with coefficients $a$ and $b$
that meet transversely at a point $p$, the pair 
$(X,\Delta)$ is klt if and only
if $a<1$ and
$b<1$. Let $\widetilde{X}$ be the blow-up of $X$
at $p$; then the coefficient in $\widetilde{\Delta}$ of the exceptional
curve $E$ is $a+b-1$. Therefore the terminal model of $(X,\Delta)$ will
blow up the point $p$ exactly when $a+b-1\geq 0$. In fact, if $a$
and $b$ are close to 1, then $a+b-1$ is also close to 1, although slightly
smaller. So the terminal model of $(X,\Delta)$ may involve arbitrarily
many blow-ups, depending on how close
the coefficients $a$ and $b$ are to 1.

We are given a terminal Calabi-Yau pair $(X,\Delta)$.
If $\Delta=0$, then $X$ is a smooth Calabi-Yau surface
and we know the cone conjecture by Sterk-Looijenga-Namikawa
\cite{Sterk, Kawamata}. So we can assume that $\Delta\neq 0$.
Using the minimal model program for surfaces, Nikulin showed
that $X$ is either rational or a $\P^1$-bundle over an elliptic
curve with $\Delta$ nef \cite[4.2.1]{NikulinMPI},
\cite[Lemma 1.4]{AM}. In the latter case, the nef effective cone
is rational polyhedral in $N^1(X)\cong \R^2$, spanned by a fiber
of the $\P^1$-bundle
together with $\Delta$ (which gives an elliptic fibration
of $X$). Thus the cone conjecture is true for $X$.

Thus, from now on, we can assume that the smooth
projective surface $X$ is rational. One consequence is that
$\Pic(X)\otimes_{\Z}\R =N^1(X)$; that is, we need not distinguish
between linear and numerical equivalence on $X$. We also deduce
that rational points in the nef cone are effective, as follows.

\begin{lemma}
\label{effectiveRR}
Let $X$ be a smooth projective rational surface with
$-nK_X$ effective for some $n>0$. Let $L$ be
a nef line bundle on $X$. Then $L$ is effective.
\end{lemma}

{\bf Proof. }Since $X$ is a rational variety, the holomorphic
Euler characteristic $\chi(X,O)$ is 1. By Riemann-Roch,
$\chi(X,L)=(L^2+L\cdot (-K_X))/2 + 1$. Since $L$ is nef
and a multiple of $-K_X$ is effective, we have $\chi(X,L)\geq 1$.
An effective divisor equivalent to $-nK_X$ is nonzero, since
$X$ is rational. Since $L$ is nef, it follows that 
an ample line bundle $A$ has $A\cdot (K_X-L)<0$.
Therefore
$h^0(X,K_X-L)=0$. Thus $h^0(X,L)=h^0(X,L)+h^0(X,K_X-L)\geq
\chi(X,L)\geq 1$. \qed

If $X$ has Picard number at most 2, then $A^e(X)$ is rational polyhedral
and so the cone conjecture is true. (For Picard number 2, since $X$
is rational, it is a $\P^1$-bundle over $\P^1$,
$X\cong P(O\oplus O(a))$ for some $a \geq 0$. The nef effective cone
is spanned by two semi-ample divisors, corresponding to the projection
$X\arrow \P^1$ and the contraction of the $(-a)$-section.)
From now on, we can assume that
$X$ has Picard number at least 3. We do this to ensure
that every $K_X$-negative extremal ray in $\NE(X)$ is spanned
by a $(-1)$-curve \cite[Lemma 1.28]{KM}.

We are assuming that $K_X+\Delta\equiv 0$,
and so $-K_X\equiv \Delta$ is effective. As a result, $-K_X$ has a Zariski
decomposition $-K_X=P+N$, meaning that $P$ is a nef $\Q$-divisor class,
$N$ is an effective $\Q$-divisor
with negative definite intersection pairing
among its components, and $P\cdot N=0$ \cite[Theorem 14.14]{Badescu}.
I claim that $P$ is semi-ample.
Indeed, the properties stated of the Zariski decomposition
imply that the effective $\R$-divisor $\Delta$ numerically equivalent
to $-K_X$ must contain $N$;
that is, the divisor $P:=\Delta-N$ is effective.
By Lemma \ref{abun}, every nef
effective $\R$-divisor on $X$ is semi-ample.
Thus $P$ is semi-ample.

In particular, the Iitaka dimension of $P$ is either 0, 1, or 2, and this
gives the main division of the proof into cases. (By definition,
$P$ has Iitaka
dimension $r$ if there is a positive integer $N$ and positive
numbers $a,b$ such that
$am^r\leq h^0(X,mP)\leq bm^r$ for all positive multiples $m$ of $N$
\cite[Chapter 10]{Iitaka}.)
The sections of multiples
$mP$ can be identified with the sections of $-mK_X$ by adding $mN$, when
$-mK_X$ and $mP$ are both
integral divisors;  so we can also describe the three cases
as $-K_X$ having Iitaka dimension 0, 1, or 2. The group
$\Aut^*(X)=\im(\Aut(X)\arrow GL(N^1(X)))$ is finite when $-K_X$
has Iitaka dimension 2 and virtually abelian for Iitaka dimension 1,
whereas it can be a fairly general group acting on hyperbolic space
when the Iitaka dimension is 0.

We start with the easiest case, where $-K_X$ is big
(that is, it has Iitaka dimension 2). In this case, we will show
that the Cox ring $\Cox(X)\cong \oplus_{L\in \Pic(X)}H^0(X,L)$
is finitely generated,
which is stronger than the cone conjecture. 

The following argument
works in any dimension. Since $\Delta$ is big, it is $\R$-linearly
equivalent to $A+E$ for an ample $\R$-divisor $A$ and an effective
$\R$-divisor $E$. Let $\Gamma=(1-\epsilon)\Delta+\epsilon E$
for $\epsilon>0$ small. Then $\Gamma$ is effective,
$(X,\Gamma)$ is klt, and $-(K_X+\Gamma)\equiv \epsilon A$
is ample. That is, $(X,\Gamma)$ is a klt Fano pair.
Birkar-Cascini-Hacon-M\textsuperscript{c}Kernan showed
that klt Fano pairs of any dimension have finitely generated
Cox ring \cite[Corollary 1.3.1]{BCHM}, as we want.

In dimension 2, this was known earlier: by the cone theorem,
a klt Fano pair $(X,\Gamma)$ has rational polyhedral cone of curves,
and every face of this cone can be contracted \cite[Theorem 3.7]{KM}.
In dimension 2, that is enough (together with $\Pic(X)\otimes\R=
N^1(X)$) to imply that the Cox ring is finitely
generated \cite{HK}. In particular, the nef effective
cone is rational polyhedral, and so the cone conjecture
is true for $X$.

Next, suppose that $-K_X$ has Iitaka dimension 1. Then the semi-ample
divisor $P$ determines a fibration of $X$ over a curve $B$, and
$P^2=0$. We have $-K_X\cdot P=(P+N)\cdot P=0$, and so the generic fiber
of $X\arrow B$ has genus 1. 
By repeatedly contracting $(-1)$-curves
contained in the fibers of $X\arrow B$, we find a factorization
$X\arrow Y\arrow B$ through a minimal elliptic surface
$Y\arrow B$. (The curves being contracted need not
be those in $N$, as one sees in examples.)
Write $\pi$ for the contraction $X\arrow Y$
and $\Delta_Y=\pi_*(\Delta)$. 
Since $K_X+\Delta\equiv 0$, we have $K_X+\Delta=\pi^*
(K_Y+\Delta_Y)$. So $(Y,\Delta_Y)$ is a klt Calabi-Yau pair
and $(X,\Delta)$ is the terminal model of $(Y,\Delta_Y)$.

We know the cone conjecture for the minimal
rational elliptic surface $Y$ \cite[Theorem 8.2]{Totaro}.
But in general, blowing up a point on a surface can increase the complexity
of the nef cone,
for example turning a finite polyhedral cone into one which is not
finite polyhedral. Rather than reduce to that earlier result,
we will go through the argument directly for $X$.

Since $-K_Y\equiv\Delta_Y$, where $K_Y$ has degree zero
on all curves contracted by $Y\arrow B$, $\Delta_Y$
is the sum of some positive multiples of fibers of $Y\arrow B$.
(Here a fiber means the pullback to $Y$ of a point in $B$,
as a divisor. We are using that the intersection pairing on
the curves contained in a fiber is negative semidefinite,
with radical spanned by the whole fiber \cite[Lemma 8.2]{BPV}.) 
The Mordell-Weil group of the elliptic fibration
$X\arrow B$ is defined as the group $\Pic^0(X_{\eta})$ where
$X_{\eta}$ is the
generic fiber. The Mordell-Weil group
acts by birational automorphisms on $Y$ over $B$,
hence by automorphisms of $Y$ since $Y\arrow B$ is minimal.
By our description of $\Delta_Y$, the Mordell-Weil group
preserves $\Delta_Y$ on $Y$. Since $(X,\Delta)$ is the
terminal model of $(Y,\Delta_Y)$ (and terminal models are
unique in dimension 2), 
the Mordell-Weil group acts by automorphisms of $(X,\Delta)$.

The main problem is to show that $\Aut(X,\Delta)$ has finitely many orbits
on the set of $(-1)$-curves in $X$. To do that, we first show that
for each curve $C$ in a fiber of $X\arrow B$, the intersection
number of a $(-1)$-curve $E$ with $C$ is bounded, independent of $E$.
Indeed, we have $1=-K_X\cdot E=(P+N)\cdot E$. So, if $E$ is not
one of the finitely many curves in $N$, we have $N\cdot E\geq 0$
and hence $P\cdot E\leq 1$. Since each fiber of $X\arrow B$ is numerically
equivalent to a multiple of $P$, this gives a bound for $E\cdot C$
for each curve $C$ contained in a fiber of $X\arrow B$, as we want.

The Picard group $\Pic (X_{\eta})$ is the quotient
of $\Pic (X)$ by some class $aP$, $a>0$, together with
all the curves $C_1,\ldots,C_r$ in reducible fibers of $X\arrow B$.
(Indeed, if a fiber contains only one curve $C$, then the class of $C$
in $\Pic(X)$ is some positive multiple of $P$.)
The degree of a line bundle on $X$
on a general fiber of $X\arrow B$ is given by the intersection
number with $bP$, for some $b>0$,
and so the Mordell-Weil group $G:=\Pic^0(X_{\eta})$
is the subquotient of $\Pic (X)$ given by
$$G=P^{\perp}/(aP, C_1,\ldots,C_r).$$

An element $x$ of the group $G$
acts by a translation on the curve $X_{\eta}$ of genus 1,
which extends to an automorphism of $X$ as we have shown.
This gives an action of $G$ on $\Pic (X)$.
We know how translation by an element $x$ of $\Pic^0(X_{\eta})$
acts on $\Pic (X_{\eta})$: by $\varphi_x(y)=
y+\deg(y)x$. Since $bP\in \Pic (X)$ is the class of a general
fiber of $X\arrow B$, this means that $\varphi_x$ acts on $\Pic(X)$
by $$\varphi_x(y)=y+(y\cdot bP) x \pmod{aP,
C_1,\ldots,C_r}.$$ 

For an element $x\in P^{\perp}$ with $x\cdot C_i=0$ for all
the curves $C_i$, the automorphism $\varphi_x$ of the minimal elliptic
surface $Y$ acts on all singular fibers $F$ of $Y\arrow B$ by
an automorphism in the identity component of $\Aut(F)$
\cite[Prop.\ 8.12(iii)]{Raynaud}. This implies the same statement
on the surface $X$ (where we have blown up some singular points
of fibers of $Y\arrow B$). In particular, $\varphi_x$ gives
the identity permutation of the curves in each fiber of $X\arrow B$.
Using that the action of $G$ on $\Pic(X)$ preserves
the intersection product, we deduce that $\varphi_x$
acts on $\Pic(X)$ by the strictly parabolic transformation
$$\varphi_x(y)=y +(y\cdot bP )x - \big[ x\cdot y
+(1/2)(x\cdot x) (y\cdot bP) \big] (bP)$$
for all $x\in P^{\perp}$ with $x\cdot C_i=0$ for all
the curves $C_i$, and for all $y\in \Pic(X)$.

Now let $E_1$ and $E_2$ be any two $(-1)$-curves on $X$
not contained in fibers of $X\arrow B$
such that $E_1\cdot bP=E_2\cdot bP$ (write $m=E_1\cdot bP$),
$E_1\cdot C_i=E_2\cdot C_i$ for all the curves
$C_i$, and $E_1\equiv E_2\pmod{m\Pic(X)}$.
Let $x=(E_2-E_1)/m\in \Pic(X)$. Then $x$ is in $P^{\perp}$,
we have
$x\cdot C_i=0$ for all the curves $C_i$, and
\begin{align*}
\varphi_x(E_1) &= E_1+(E_1\cdot bP) x
- \big[ (x\cdot E_1) + (1/2)(x\cdot x)(E_1\cdot bP) \big] (bP)\\
 &= E_2.
\end{align*}

We have shown that the intersection numbers $m=E\cdot bP$
and $E\cdot C_i$ are bounded, among all $(-1)$-curves $E$
on $X$.
So, apart from the finitely many $(-1)$-curves
contained in fibers of $X\arrow B$,
the $(-1)$-curves $E$ are divided into finitely many
classes according to $m$, the intersection numbers of $E$
with the curves $C_i$,
and the class of $E$ in $\Pic(X)/m$. By the previous paragraph,
the $(-1)$-curves on $X$ fall into finitely many
orbits under the action of $G$ we defined. This completes
the proof that $\Aut(X,\Delta)$ has finitely many orbits
on the set of $(-1)$-curves.

We now describe all the extremal rays of the cone of curves $\NE(X)$,
following Nikulin \cite[Proposition 3.1]{NikulinIJM}.
We have mentioned that every $K_X$-negative extremal ray
is spanned by a $(-1)$-curve. On the other side,
a $K_X$-positive extremal ray must be spanned by
one of the finitely many curves in $N$. (Since $P\cdot N=0$,
the curves in $N$ are contained in fibers of the elliptic
fibration $X\arrow B$ given by $P$.) Finally, let
$\R^{\geq 0}x$ be an extremal ray of $\NE(X)$ in $K_X^{\perp}$.
Suppose $x$ is not a multiple of a curve in $N$; then
$N\cdot x\geq 0$ and $P\cdot x=0$. Since $-K_X\equiv P+N$,
it follows that $P\cdot x=0$. Since $P^2=0$, the Hodge index theorem
gives that $x$ is a multiple of $P$ in $N^1(X)$ or $x^2<0$.
In the latter case, the ray $\R^{\geq 0}x$ must be spanned by
a curve $C$. Since $P\cdot C=0$, $C$ is a curve in some fiber of
$X\arrow B$. There are only finitely many numerical equivalence
classes of curves in the fibers of $X\arrow B$. We conclude
that almost all (all but finitely many) extremal rays of $\NE(X)$
are spanned by $(-1)$-curves.

Moreover, the only possible limit ray of the $(-1)$-rays is
$\R^{\geq 0}P$. Indeed, if $\R^{\geq 0}x$ is a limit ray of $(-1)$-rays,
then $x^2=0$ and $-K_X\cdot x=0$. (For an ample line bundle $A$,
all $(-1)$-curves have $E^2=-1$ and $-K_X\cdot E=-1$, while
their degrees $A\cdot E$ in an infinite sequence must tend to infinity.
Since $0<A\cdot x<\infty$, this proves the properties stated of $x$.)
Also, $N\cdot x\geq 0$ since $N$ has nonnegative intersection
with almost all $(-1)$-curves, and $P\cdot x\geq 0$, while
$-K_X=P+N$; so $P\cdot x=0$. Since $P^2=0$ and $x^2=0$,
$x$ is a multiple of $P$ by the Hodge index theorem.

We can deduce that the nef cone $\A(X)$ is rational polyhedral
near any point $y$ in $\A(X)$ not in the ray $\R^{>0}P$. First, such 
a point has $y^2\geq 0$ and also $P\cdot y\geq 0$, since
$y$ and $P$ are nef. If $P\cdot y$ were zero, these properties
would imply that $y$ was a multiple of $P$; so we must
have $P\cdot y>0$. Since the only possible limit ray of $(-1)$-rays
is $\R^{\geq 0}P$, there is a neighborhood of $y$ which has
positive intersection with almost all $(-1)$-curves. 
Since almost all extremal rays of $\NE(X)$ are spanned
by $(-1)$-curves, we conclude that the nef cone $\A(X)$
is rational polyhedral near $y$, as claimed.

In particular, for each $(-1)$-curve $E$ not contained
in a fiber of $X\arrow B$, the face $\A(X)\cap E^{\perp}$
of the nef cone is rational polyhedral, since it does not
contain $P$. So the cone $\Pi_E$ spanned by $P$
and $\A(X)\cap E^{\perp}$ is rational polyhedral.

\begin{floatingfigure}[l]{0.3\textwidth}
\centering
\includegraphics[scale=0.3]{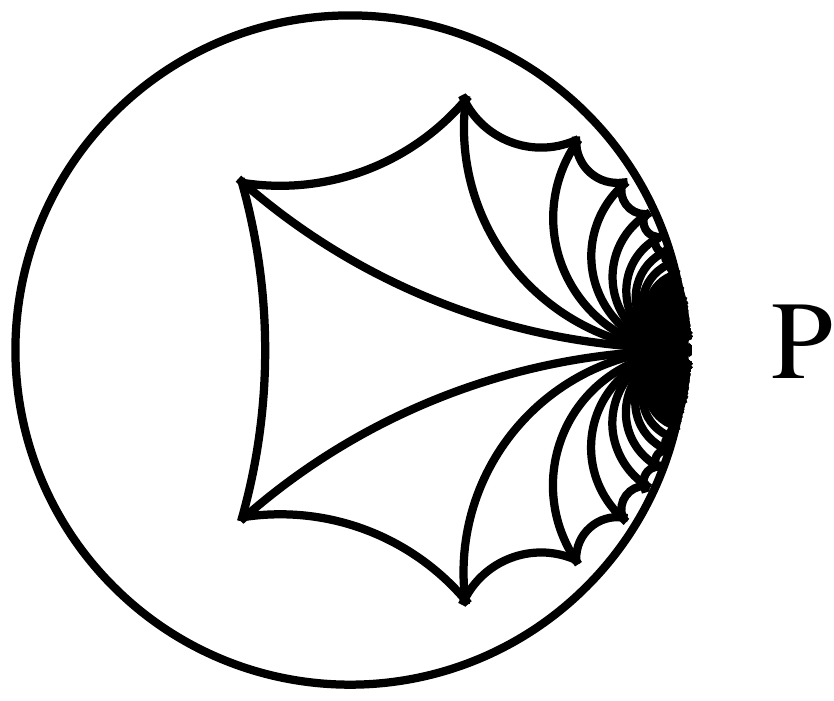}
\end{floatingfigure}
Let $x$ be any nef $\R$-divisor on $X$.
Let $c$ be the maximum real number such that $y:=x-cP$
is nef. Then $x$ and $y$ have the same degree on all curves
contained in a fiber
of $X\arrow B$. By our list of the extremal
rays in $\NE(X)$,
there must be some $(-1)$-curve $E$ not contained in a fiber
such that $y\in E^{\perp}$.
Therefore $x$ is in the cone $\Pi_E$. That is, the nef cone
$\A(X)$ is the union of the rational polyhedral cones
$\Pi_E$, as in the figure.

Any rational point $x$ in the nef cone $\A(X)$ is effective,
by Lemma \ref{effectiveRR}. So the rational polyhedral cones
$\Pi_E$ are contained in $A^e(X)$, and $A^e(X)=\A(X)$
is the union of these cones. Since there are
only finitely many $\Aut(X,\Delta)$-orbits
of $(-1)$-curves $E$, Lemma \ref{dirichlet} proves
the cone conjecture for $X$.

It remains to consider the case where
 $-K_X$ has Iitaka dimension 0. In the Zariski decomposition
$-K_X\equiv\Delta=P+N$, $P$ is numerically trivial, so
$-K_X\equiv N$ where $N$ is an effective $\Q$-divisor with negative
definite intersection pairing on its irreducible components.
In this case, $N$ is the unique effective
$\R$-divisor numerically equivalent to $-K_X$, and so the given
divisor $\Delta$ is equal to $N$.

We use the following negativity lemma, which is essentially
an elementary result on quadratic forms \cite[section V.3.5]{Bourbaki}.

\begin{lemma}
\label{neg}
(Negativity lemma) Let $N$ be a set of curves
on a smooth projective surface on which the intersection
pairing is negative definite. Let $-D$ be a linear combination
of the curves $N_i$ which has nonnegative intersection with each $N_i$.
Then $D$ is effective. Moreover, the support of $D$ is a union
of some connected components of $N$.
\end{lemma}

In our case, we can contract all the curves in $N$.
(There is a positive linear combination
$D=\sum a_iN_i$ with $D\cdot N_i=-1$ for all $i$ by Lemma \ref{neg}.
The pair $(X,\Delta+\epsilon D)$ is klt for
$\epsilon>0$ small, and $(K_X+\Delta+\epsilon D)\cdot N_i=-\epsilon <0$
for all $i$, so we can contract the $N_i$'s by the cone theorem
\cite[Theorem 3.7]{KM}.) Write $\pi:X\arrow Y$ for the resulting
contraction. 
Since $K_X+\Delta\equiv 0$, we have $K_X+\Delta=\pi^*(K_Y)$.
So $Y$ is a klt Calabi-Yau surface
and $(X,\Delta)$ is the
terminal model of $Y$.

We know the cone conjecture
for $Y$ by Theorem \ref{surface}.
But that does not immediately imply the
statement for $X$. In general, blowing up a point on a surface increases
the Picard number and can
make the nef cone more complicated, for example turning a finite polyhedral
cone into one which is not finite polyhedral.
Since $(X,\Delta)$ is the
terminal model of $Y$ (and terminal models are
unique in dimension 2), every automorphism of $Y$
lifts to an automorphism of $(X,\Delta)$.
Thus it will suffice to show that $\Aut(Y)=\Aut(X,\Delta)$
has a rational polyhedral fundamental domain
on the nef effective cone of $X$.

We can describe all the extremal rays of the cone of curves $\NE(X)$,
following Nikulin \cite[Proposition 3.1]{NikulinIJM}.
We have mentioned that every $K_X$-negative extremal ray
is spanned by a $(-1)$-curve. On the other side,
a $K_X$-positive extremal ray $\R^{\geq 0}x$ must be spanned by
one of the finitely many curves in $N$, since
$0>-K_X\cdot x= N\cdot x$. Finally, let
$\R^{\geq 0}x$ be an extremal ray of $\NE(X)$ in $K_X^{\perp}$.
This ray may be spanned by one of the curves $N_i$. If it is not,
then $x\cdot N_i\geq 0$ for all $i$. Therefore $0=-K_X\cdot x
=N\cdot x\geq 0$, and so $N_i\cdot x=0$ for all $i$.
That is, $x=\pi^*(w)$ for some $w\in \NE(Y)$.

Let $N_1,\ldots,N_r$ be the irreducible
components of $N$. A $(-1)$-curve $C$
in $X$ has $1=-K_X\cdot C= (\sum a_iN_i)\cdot C=\sum a_i\lambda_i$,
where $a_1,\ldots,a_r$ are fixed positive numbers and we
write $\lambda_i=C\cdot N_i$. As a result,
there are only finitely many possibilities
for the natural numbers $(\lambda_1,\ldots,\lambda_r)$,
for all $(-1)$-curves on $X$ not among the curves $N_i$.
Call these the finitely many
{\it types }of $(-1)$-curves on $X$.

We now describe the nef cone of $X$.
Every divisor class $u$ on $X$ can be written as $\pi^*(y)-\sum b_iN_i$
for some real numbers $b_i$ and some $y\in N^1(Y)$. If $u$ is nef,
then $y$ must be nef on $Y$. Also, $u$
has nonnegative degree on the curves $N_i$, which says that
$(b_1,\ldots,b_r)$ lies in a certain rational polyhedral cone
$B$. The cone $B$ is contained in $[0,\infty)^r$
by the negativity lemma, Lemma \ref{neg}.
By our description of the extremal rays of $\NE(X)$,
a class $u=\pi^*(y)-\sum b_iN_i$ in $N^1(X)$
is nef if and only if $y$ is nef
on $Y$, $(b_1,\ldots,b_r)$ is in the cone $B$, and $u$ has nonnegative degree
on all $(-1)$-curves not among the curves $N_i$ in $X$.
The last condition says, more explicitly: for each $(-1)$-curve 
$C$ not among the curves $N_i$ in $X$, we must have
\begin{align*}
0 &\leq C\cdot [\pi^*(y)-\sum b_iN_i]\\
&= y\cdot \pi_*(C)-\sum \lambda_ib_i
\end{align*}
where we write $\lambda_i=C\cdot N_i$. 

Thus, for $u=\pi^*(y)-\sum b_iN_i$ to be nef means that the
numbers $b_i$ satisfy the upper bounds $\sum \lambda_ib_i\leq
y\cdot \pi_*(C)$ for all $(-1)$-curves $C$ on $X$ not among
the curves $N_i$, where $(\lambda_1,\ldots,\lambda_r)$ is
the type of $C$.
Notice that a $(-1)$-curve $C$
on $X$ is determined by its type together with the class $\pi_*(C)$ in
$N_1(Y)$.

Since $\Aut(X,\Delta)=\Aut(Y)$,
the theorem holds if for every rational polyhedral cone
$S\subset A^e(Y)$, the inverse image $T$ of $S$ under
$\pi_*: A^e(X)\arrow A^e(Y)$
is rational polyhedral. (Then the inverse images of any decomposition
given by the cone conjecture for $Y$ form a decomposition satisfying
the cone conjecture for $(X,\Delta)$.)
Let us first define $T$ to be the inverse
image of $S$ in the nef cone $\A(X)$; at the end we will
check that $T$ is actually contained in the nef effective cone.
Since $S-0$ is compact modulo scalars, it suffices
to prove that $T$ is rational polyhedral in
the inverse image of some neighborhood of each
nonzero point in $S$. 

First, let $y_0\in S$ be a point with
$y_0^2>0$. We want to show that only finitely many
$(-1)$-curves in $X$ are needed to define the cone $T$ over a neighborhood
of $y_0$ in $S$. It suffices to show that for each type
$\lambda$ of $(-1)$-curves on $X$, there is a finite set $Q$
of $(-1)$-curves of type $\lambda$ such that for all
$y$ in some neighborhood of $y_0$, $y\cdot \pi_*(C)$ is minimized
among all $(-1)$-curves $C$ of type $\lambda$ by one of the
curves in $Q$. The point is that the type of the $(-1)$-curve
determines the rational number $c:=\pi_*(C)^2$. (This can be positive,
negative, or zero, as examples show.) The intersection pairing
on $N_1(Y)$ has signature $(1,*)$ by the Hodge index theorem.
Since $y_0^2>0$, the intersection of the hyperboloid
$\{ z\in N_1(Y): z^2=c\}$ with $\{z\in N_1(Y): |z\cdot y_0|\leq M\}$
is compact, for any number $M$.
So there are only finitely many integral classes $z$ in $N_1(Y)$
with $z^2=c$ and with given bounds on $z\cdot y_0$, and
the same finiteness applies for $y$ in some neighborhood of $y_0$.
Thus only finitely many classes $z=\pi_*(C)$, and hence
only finitely many $(-1)$-curves $C$, can minimize
$y\cdot \pi_*(C)$ for any $y$ in a neighborhood of $y_0$, as we want.

It remains to consider a nonzero 
point $y_0$ in $S$ with $y_0^2=0$. Since $S\subset
A^e(Y)$ is a rational polyhedral cone contained in the positive
cone $\{y\in N^1(Y): y^2\geq 0, A\cdot y\geq 0\}$, $y_0$ must
belong to an extremal ray of $S$. Therefore we can scale $y_0$
to make it an integral
point in $N^1(Y)$ (the class of a line bundle on $Y$).
Since $y_0$ is a nef integral divisor on
the klt Calabi-Yau surface $Y$, it is semi-ample by
Lemmas \ref{effectiveRR} and \ref{abun}. Since $y_0^2=0$, the
corresponding contraction maps $Y$ onto a curve $L$.

For each point $p$ of $Y$ over which $\pi:X\arrow Y$ is not
an isomorphism, let $D$ be a curve through $p$
which is contained in a fiber of $Y\arrow L$ (clearly there is
such a curve). Let $C$ be the proper transform of $D$ in $X$.
Since $C$ is contained in a singular fiber of $X\arrow L$,
$C$ has negative self-intersection and hence spans an
extremal ray of $\NE(X)$. Since $C$ is not among the curves $N_i$
and has positive intersection with some curve $N_i$, our
description of the extremal rays of $\NE(X)$ shows that
$C$ is a $(-1)$-curve. Thus, for each connected component $R$ of $N$
(corresponding to a point over which $\pi:X\arrow Y$
is not an isomorphism), there is a $(-1)$-curve
$C$ on $X$ such that $y_0\cdot \pi_*(C)=0$ and $\lambda_i=C\cdot N_i$
is positive for some $N_i$ in $R$.

Moreover, the set $Q$ of $(-1)$-curves $C$ in $X$ with $y_0\cdot
\pi_*(C)=0$ is finite, since such a curve must be contained
in one of the
finitely many singular fibers of $X\arrow L$. 
I claim that these finitely many $(-1)$-curves are enough
to define the cone $T$ over a neighborhood of the vertex $y_0$ in
the rational polyhedral cone 
$S$. We can view such a neighborhood (up to scalars)
as the set of linear combinations $y=y_0+\sum c_i v_i$,
for some nef classes $v_i$ on $Y$, with $c_i\geq 0$ near zero.
Therefore $y\cdot \pi_*(C)\geq y_0\cdot \pi_*(C)$ for
all $(-1)$-curves  $C$ in $X$. So $y\cdot \pi_*(C)$ is at least 1
for the $(-1)$-curves $C$ outside the set $Q$, whereas it
is small (for $c_i$ near zero) for $C$ in the set $Q$.
Therefore the inequality $\sum \lambda_i b_i\leq 
y\cdot \pi_*(C)$
is only needed for the finitely many curves $C$ in $Q$; that is,
$T$ is rational polyhedral over a neighborhood of $y_0$ in $S$.

To check this in detail, we have to recall our earlier
comment that for each connected component
$R$ of $N$, $Q$ contains
a $(-1)$-curve $C$ with $\lambda_i>0$ for some $N_i$ in $R$.
This is needed to show
that the inequalities for $C$ in $Q$ imply the inequalities
for all $(-1)$-curves $C$ in $X$. Namely, the inequalities for
$C$ in $Q$ imply that $b_i$ is small (assuming $y$ is near $y_0$)
for some $N_i$ in each
connected component of $N$. By the negativity
lemma, since $(b_1,\ldots,b_r)$ is in $B$, it follows that every $b_i$
is small. Indeed, Lemma \ref{neg}
says that if a point
in the rational polyhedral cone $B$ has one $b_i$ equal to zero,
then $b_j$ is also zero for every $N_j$ in the same connected
component as $N_i$.
This implies the same statement for ``small''
in place of ``zero''.

Thus the cone $T\subset \A(X)$ is rational polyhedral. We actually
want to know that this rational polyhedral cone is contained
in $A^e(X)$. That is the case, by Lemma \ref{effectiveRR}
(on a smooth projective rational surface, every nef $\Q$-divisor
class is effective). As explained earlier, since 
$T\subset A^e(X)$ is rational polyhedral,
the cone conjecture for $(X,\Delta)$ is proved.
\qed

\section{Finite generation of the Cox ring}

\begin{corollary}
\label{finite}
Let $(X,\Delta)$ be a klt Calabi-Yau pair of dimension 2
over the complex numbers. The following are equivalent:

(1) The nef effective cone of $X$ is rational polyhedral.
(This means in particular that the nef effective cone is closed.)

(2) The nef cone of $X$ is rational polyhedral.

(3) The image of $\Aut(X)\arrow GL(N^1(X))$ is a finite group.

(4) The image of $\Aut(X,\Delta)\arrow GL(N^1(X))$ is a finite group.

If the first Betti number of $X$ is zero, properties (1) to (4)
are equivalent to finite generation of the Cox ring
$\Cox(X)\cong \oplus_{L\in \Pic(X)}H^0(X,L)$.
\end{corollary}

For this class of varieties, property (4) is often
an easy way to determine whether the Cox
ring is finitely generated. For example, for minimal
rational elliptic surfaces,
property (4) is equivalent to finiteness of the Mordell-Weil
group, which can be described in simple geometric terms
\cite[Theorem 5.2, Theorem 8.2]{Totaro}. The
rational elliptic surfaces with finite
Mordell-Weil group have been classified
by Miranda-Persson \cite{MP} and Cossec-Dolgachev
\cite{CD}. (See Prendergast-Smith \cite{PS} for an analogous
classification in dimension three.) The K3 surfaces and Enriques
surfaces with finite automorphism group were
classified by Nikulin, Vinberg, and Kondo \cite{NikulinICM, 
Nikulincone, Kondo}.

In any dimension, every variety $X$ of {\it Fano type}, meaning
that there is a divisor $\Gamma$ with $(X,\Gamma)$ klt Fano, has
finitely generated Cox ring,
by Birkar-Cascini-Hacon-M\textsuperscript{c}Kernan
\cite[Corollary 1.3.1]{BCHM}. The varieties of Fano type form a subclass
of the varieties $X$ which have a divisor $\Delta$
with $(X,\Delta)$ klt Calabi-Yau, namely the subclass
with $-K_X$ big. (Compare the proof of Theorem \ref{rational}
in the case where $-K_X$ is big.) For example,
the blow-up $X$ of $\P^2$
at any number of points on a conic
is of Fano type.
Therefore $X$ has finitely generated Cox
ring, as Galindo-Monserrat \cite[Corollary 3.3]{GalM},
Castravet-Tevelev \cite{CT} and Mukai \cite{Mukaifin}
proved by other methods. Recently Testa, V\'arilly-Alvarado and Velasco
showed that
every smooth rational surface with $-K_X$ big has finitely
generated Cox ring \cite {TVV}; Chenyang Xu showed that such a surface
need not be of Fano type.

{\bf Proof of Corollary \ref{finite}. }The closure
of the nef effective cone is the nef cone,
and so (1) implies (2). The subgroup of $GL(n,\Z)$ preserving a
rational polyhedral cone that does not contain a line and has
nonempty interior is finite, and so (2) implies (3). Clearly
(3) implies (4). Theorem
\ref{rational} shows that $\Aut(X,\Delta)$ acts with rational
polyhedral fundamental domain on the nef effective cone,
and so (4) implies (1).

Suppose that these equivalent conditions hold.
Since $(X,\Delta)$ is a klt Calabi-Yau pair, every
nef effective $\Q$-divisor on $X$ is semi-ample by Lemma
\ref{abun}. Another way to say this is that every face
of the cone of curves can be contracted. If $b_1(X)=0$,
Hu and Keel showed that (since $X$ has dimension 2)
finite generation of the Cox ring is equivalent to
the nef cone being rational polyhedral together with every
nef divisor being semi-ample \cite{HK}. \qed

\small \sc DPMMS, Wilberforce Road,
Cambridge CB3 0WB, England

b.totaro@dpmms.cam.ac.uk
\end{document}